\documentclass[10pt]{article}
\usepackage[T1]{fontenc}
\usepackage[utf8]{inputenc}
\usepackage{lmodern}
\usepackage[margin=0.9in]{geometry}
\usepackage{amsmath,amssymb}
\usepackage{array,booktabs,longtable}
\usepackage{xcolor}
\usepackage{url}
\usepackage[colorlinks=true,linkcolor=blue!50!black,citecolor=blue!50!black,urlcolor=blue!50!black]{hyperref}

\DeclareUrlCommand\lean{\urlstyle{tt}}
\makeatletter
\g@addto@macro{\UrlBreaks}{\do\.\do\_}
\makeatother

\newcommand{\C}{\mathbb{C}}
\newcommand{\Z}{\mathbb{Z}}
\newcommand{\sltwo}{\mathfrak{sl}_2}
\newcommand{\slthree}{\mathfrak{sl}_3}
\newcommand{\g}{\mathfrak{g}}
\newcommand{\Vn}{V^\natural}

\title{Vertex Operator Algebras in Lean:\\ a formal library and its companion book}
\author{Fedor K.\ Popov\textsuperscript{a} \quad Leonardo Rastelli\textsuperscript{a,b}\\[0.8em]
  \small\textsuperscript{a}Simons Center for Geometry and Physics, Stony Brook University, Stony Brook, NY 11794-3840, USA\\
  \small\textsuperscript{b}C.\,N.\ Yang Institute for Theoretical Physics, Stony Brook University, Stony Brook, NY
  11794-3840, USA\\[0.5em]
  \small\texttt{fpopov@scgp.stonybrook.edu}, \texttt{leonardo.rastelli@stonybrook.edu}}
\date{Release v1.1}

\begin{document}
\maketitle

\begin{abstract}
  We present \texttt{LeanVOA}, a Lean~4 library for the theory of vertex operator algebras (VOAs), together with a
  companion book that sets each formalized result next to its textbook statement. The library constructs the
  Moonshine module $V^\natural$ and proves, with no hypotheses, the Conway--Norton moonshine conjecture for every
  automorphism of $V^\natural$ (Borcherds' theorem), along Borcherds' route through the Monster Lie algebra. It
  identifies $\mathrm{Aut}\,V^\natural$ with a formal model of the Monster taken from the AxiomMath library for the
  classification of finite simple groups. A separate package in the repository imports that library's generation
  theorem for the model and proves, with no hypotheses, that $\mathrm{Aut}\,V^\natural$ is isomorphic to the model,
  compatibly with the Griess algebras, and that moonshine holds for every element of the Monster. It also formalizes Zhu's modular invariance, Huang's tensor
  categories with the Verlinde formula, the Gannon--Rayhaun classification of $c=1$ VOAs, and classical freeness of the boundary admissible
  $L_k(\sltwo)$ (Arakawa--Dai) and of the minimal models $M(2,2m+3)$ (van Ekeren--Heluani). It also contains new
  results: graded first-order rigidity of all Virasoro minimal models and of $L_k(\sltwo)$ at every admissible level; classical freeness of $L_{-3/2}(\slthree)$, and of the other members of the Deligne series under stated
  assumptions; and a proof that graded
  unitarity of the boundary admissible $L_k(\sltwo)$ implies that of the minimal models $M(2,2n+1)$. The library
  declares no axiom of its own. This note announces release v1.1, which moves the whole library to Lean~v4.35.
\end{abstract}

\section{Introduction}

The library \texttt{LeanVOA} formalizes, in Lean~4 \cite{lean4} on top of Mathlib \cite{mathlib}, the algebraic
theory of vertex operator algebras: formal calculus and locality from first principles, the VOA axioms
and the Borcherds identity, the standard constructions (Heisenberg, Virasoro, affine, lattice), the module theory,
Zhu's algebra and modular invariance, Huang's tensor categories, the moonshine module and monstrous moonshine,
and several recent theorems from the literature (listed in Section~\ref{sec:scope}). The companion book, \emph{Vertex Operator Algebras in Lean}, is the
reviewed account of the library: for each definition and theorem it gives the mathematical statement, the Lean
declaration copied verbatim from the source, and a ``dictionary'' that decodes every convention and hypothesis
and says whether the Lean statement is the same as, stronger than, or weaker than the textbook one.

\begin{itemize}
  \item Repository: \url{https://github.com/Fedor-Popov/Lean_VOA}
  \item Release: v1.1 (tag \texttt{v1.1}; previous releases v1.06, v1.05 and v1.0); Lean~v4.35.0-rc3 with Mathlib
    commit \texttt{b63f6e8a}
  \item Monster package: the folder \texttt{MoonshineCFSG/} of the repository, a separate Lake package that depends on
    the library and on the AxiomMath library \cite{AxiomCFSG}; the library itself depends only on Mathlib
  \item Book (PDF): \texttt{Notes/book.pdf} in the repository (2{,}807 pages); its \LaTeX{} sources are in
    \texttt{Notes/Textbook/}
  \item Short version of the book (PDF): \texttt{Notes/brief.pdf} (258 pages): the main statements of each chapter,
    each with its verbatim Lean declaration and its status, without proofs or discussion, numbered as in the book
  \item License: Apache-2.0.
\end{itemize}

The library defines VOAs in the general sense of Frenkel--Lepowsky--Meurman and Lepowsky--Li \cite{FLM,LepowskyLi}
(\lean{GradedVertexOperatorAlgebra}), but develops the theory for VOAs of CFT type ($V_0=\C\mathbf 1$, $V_n=0$ for
$n<0$; \lean{VertexOperatorAlgebra}). Every theorem below about VOAs is about VOAs of CFT type; extending the theory
to the general notion is future work.

\paragraph{New in v1.05: monstrous moonshine.} Release v1.0 constructed $\Vn$ with no Monster action. Release v1.05
adds monstrous moonshine, in the library and in a new Part~VIII of the book (seven chapters). The theorem
\lean{VOA.Moonshine.ConwayNorton.conwayNorton}\,:\,\lean{ConwayNortonConjecture} has no hypotheses: for every
automorphism $g$ of $\Vn$, the McKay--Thompson series
\begin{equation*}
  T_g(\tau)=q^{-1}+\sum_{n\ge1}\mathrm{tr}\bigl(g\,|\,\Vn_{n+1}\bigr)\,q^n
\end{equation*}
is the normalised Hauptmodul of a genus-zero group $\Gamma_g$ that contains $\Gamma_0(N)$ for some $N\ge1$ and is
commensurable with $SL_2(\Z)$. This is the qualitative, genus-zero form of the conjecture of Conway and Norton
\cite{ConwayNorton}, proved by Borcherds \cite{Borcherds}: it does not identify $\Gamma_g$ with the group in Conway
and Norton's table, nor relate $N$ to the order of $g$, and ``genus zero'' is read off on paper from the Lean
definition of a normalised Hauptmodul. The proof follows Borcherds: the vertex algebra of the hyperbolic plane
$I\!I_{1,1}$, the tensor theory $\Vn\otimes V_{I\!I_{1,1}}$ of central charge $26$, the no-ghost theorem
\cite{GoddardThorn}, the Monster Lie algebra as a Borcherds--Kac--Moody-type Lie algebra with a free positive part,
the denominator identity (as the combination of two Lean theorems) and the twisted denominator identities, and
complete replicability of the $T_g$. Borcherds' last step, a class-by-class comparison that needs the character table
of the Monster, is replaced by a case-free proof that replicable families are Hauptmoduln (genus zero on an abstract
interface), in the spirit of Cummins and Gannon \cite{CumminsGannon}, followed by an upgrade from $\Gamma_1(N)$ to
$\Gamma_0(N)$. Along the way the library proves, again with no hypotheses and without group theory, that
$\mathrm{Aut}\,\Vn$ is finite (it acts faithfully on the finitely many rigid idempotents of the Griess algebra $\Vn_2$), that every automorphism has
finite order, and that every graded trace is rational.

\paragraph{The Monster.} For the Monster we use the rational model of the AxiomMath Lean library for the
classification of finite simple groups \cite{AxiomCFSG}, whose definitions we port: \lean{AutO} is the group of
automorphisms of their rational Griess algebra \lean{SpaceB}, of dimension $196884$, built after Conway
\cite{Conway85} from Griess' algebra \cite{Griess}. The library proves, with no hypotheses, that the Griess algebra
$\Vn_2$ is $\C\otimes_{\mathbb Q}\mathrm{SpaceB}$, with its product and form (\lean{griessIso}); that every element of AxiomMath's generated
group $\mathrm{Gmat}=\langle G_{\mathrm{vis}},\tau_B\rangle$ (with $G_{\mathrm{vis}}$ the visible group
$2^{1+24}.\mathrm{Co}_1$, classical and not formalised as such) lifts to an automorphism of $\Vn$ and satisfies
moonshine; and that $\mathrm{Aut}\,\Vn$ is isomorphic to the subgroup of \lean{AutO} that acts on $\Vn$. The last
statement is proved without Tits' theorem \cite{Tits}, from rationality of traces, Schur's lemma, and a commutant
lemma: every $\mathbb Q$-linear map commuting with $\mathrm{Gmat}$ is $a\cdot1+b\,P_0$, where $P_0$ is the projection onto
the line of the Griess unit, with the transitivity inputs from the Mathieu group $M_{24}$ checked by the kernel. Finally, given
the named hypothesis \lean{GmatEqAutO}, defined as the proposition $\mathrm{Gmat}=\mathrm{AutO}$, the theorem
\lean{monstrousMoonshine} provides an isomorphism $\varphi:\mathrm{AutO}\simeq\mathrm{Aut}\,\Vn$ compatible with the
Griess algebras, and moonshine for every element of the Monster. The hypothesis is a theorem of the AxiomMath library,
which the library does not import (Section~\ref{sec:trust}).

\paragraph{The hypothesis discharged.} The package \texttt{MoonshineCFSG/} of the repository depends on the library
and on AxiomMath's library, unmodified, at commit \texttt{5437cdc313}. It contains a bridge module and copies of the 85
library modules the bridge needs, in which the Monster definitions are built on AxiomMath's own structures, so that
\lean{Generated.Gmat} and \lean{AutO} there are definitionally AxiomMath's. AxiomMath's theorem
\lean{CFSG.Blueprint.NativeSporadic.MonsterModel.gmat_eq_griess} then proves \lean{GmatEqAutO}, and the theorem
\lean{VOA.Moonshine.ConwayNorton.monstrousMoonshine_unconditional} states monstrous moonshine for AxiomMath's
Monster with no hypotheses:
\begin{gather*}
  \exists\,\varphi:\mathrm{AutO}\simeq\mathrm{Aut}\,\Vn\quad\text{with}\quad
  \varphi(g)|_{\Vn_2}=\C\otimes g\ \text{ for all } g\in\mathrm{AutO},\\
  \text{and } T_{\varphi(g)} \text{ a normalised genus-zero Hauptmodul for every } g\in\mathrm{AutO}.
\end{gather*}
Its axioms are Lean's three standard ones (\lean{propext}, \lean{Classical.choice}, \lean{Quot.sound}). The
generation theorem $\mathrm{Gmat}=\mathrm{AutO}$ is AxiomMath's \cite{AxiomCFSG}; our contribution is the bridge that
makes it apply to the moonshine theorems. Release v1.06 first proved this theorem, in a Lean~v4.35 companion folder
next to a Lean~v4.30 library.

\paragraph{New in v1.1: one toolchain.} Release v1.1 moves the whole library to Lean~v4.35.0-rc3, the toolchain of
AxiomMath's library, with Mathlib commit \texttt{b63f6e8a}. The library depends only on Mathlib: users who do not need
the unconditional Monster theorem never download or build AxiomMath's library, which has about $200{,}000$ files and
builds about $20{,}000$ modules for the generation theorem alone. The package \texttt{MoonshineCFSG/} replaces the
companion folder of v1.06. Release v1.1 also adds a vacuity audit, prompted by an external reviewer's report: a
registry of the types the library proves empty, a ``vacuous --- do not cite'' flag on the 97 declarations whose
hypotheses mention one (none is a headline), and, for each of the six headline theorems on $C_2$-cofinite VOAs, an
instance at a concrete VOA with every hypothesis discharged. Six statements that earlier editions of the book
attributed to the literature only are now Lean theorems.

\section{Scope}\label{sec:scope}

Table~\ref{tab:scope} lists the headline results: classical theorems, recent theorems from the literature, and
new results. Each entry says what is proved and what, if anything, it assumes; the book gives the precise statements,
conventions and caveats. Every listed declaration depends only on Lean's three standard axioms (\lean{propext},
\lean{Classical.choice}, \lean{Quot.sound}). One entry, the last of the moonshine block, carries an explicit named
hypothesis, discussed in Section~\ref{sec:trust}.

{\footnotesize
  \setlength{\tabcolsep}{4pt}
  \renewcommand{\arraystretch}{1.15}
  \begin{longtable}{>{\raggedright\arraybackslash}p{0.625\textwidth}>{\raggedright\arraybackslash}p{0.33\textwidth}}
    \caption{Headline results}\label{tab:scope}\\
    \toprule
    Statement and source & Lean name(s)\\
    \midrule
    \endfirsthead
    \toprule
    Statement and source & Lean name(s)\\
    \midrule
    \endhead
    \bottomrule
    \endfoot
    \multicolumn{2}{l}{\textbf{Tier 1: foundations and classical theorems}}\\
    \midrule
    \textbf{Foundations.} VOAs of CFT type, with the Borcherds identity derived from locality rather than assumed.
    Constructed examples: the Heisenberg VOA, the universal Virasoro VOA at every central charge, the universal affine
    VOA $V^k(\sltwo)$, and the lattice VOA of every positive-definite even lattice. Source: \cite{FLM,LepowskyLi}.
    & \lean{VOA.VertexOperatorAlgebra.borcherds}, \lean{VOA.heisenbergVOA}, \lean{VOA.VirModule.virasoroVOA},
    \lean{VOA.AffineVOA.affineSl2VOA_closed}, \lean{VOA.LatGen.latticeVOAGen}\\
    \addlinespace
    \textbf{Moonshine module.} The Frenkel--Lepowsky--Meurman module $V^\natural$, with no hypotheses: a simple VOA
    with $c=24$, $V^\natural_1=0$ and $\dim V^\natural_2=196884$, whose graded dimensions are the coefficients of
    $j(\tau)-744$. Source: \cite{FLM,DGM}.
    & \lean{VOA.Moonshine.moonshine}, \lean{VOA.Moonshine.hasSum_JShift}\\
    \addlinespace
    \textbf{Zhu's modular invariance.} For a $C_2$-cofinite rational VOA, the torus one-point functions of the
    irreducible modules transform according to a representation of $SL_2(\Z)$. Instances with every hypothesis checked: all
    rank-one lattice VOAs and their $\Z_2$-orbifolds. Source: \cite{Zhu,DLMtwisted,DLMalg}.
    & \lean{VOA.zhuModularity_canonical_of_C2Cofinite_admissiblyRational},
    \lean{VOA.Classification.COne.DLMRat.admissiblyRational_of_dlmRational},
    \lean{VOA.Classification.COne.zhuModularity_V2m}\\
    \addlinespace
    \textbf{Affine VOAs for every simple Lie algebra.} The universal affine VOA $V^k(\g)$ and its simple quotient
    $L_k(\g)$ at every non-critical level, with concrete models of $\mathfrak{sl}_n$ and of the exceptional Lie
    algebras. Source: \cite{Kac}.
    & \lean{VOA.AffineGeneral.universalAffineVOA}, \lean{VOA.AffineGeneral.simpleAffineQuotientVOA_isSimple}\\
    \midrule
    \multicolumn{2}{l}{\textbf{Tier 2: recent theorems}}\\
    \midrule
    \textbf{Classical freeness of $L_k(\sltwo)$ (Arakawa--Dai).} At the boundary admissible levels $k=-2+2/(2m+1)$, $m\ge1$,
    $L_k(\sltwo)$ is classically free, with $C_2$ algebra $\C[e,h,f]/(eQ^m,hQ^m,fQ^m)$, $Q$ the Casimir; also its vacuum
    character, in product and fermionic form. Source: \cite{ArakawaDai}.
    & \lean{VOA.AffineSl2BoundaryClassicalFreeness.arakawaDai_theorem_1_1},
    \lean{VOA.AffineSl2BoundaryCharacter.character_eq_fermionic_chain}\\
    \addlinespace
    \textbf{Classical freeness of the minimal models $M(2,2m+3)$.} For every $m\ge0$, with $C_2$ algebra
    $\C[T]/(T^{m+1})$. Source: \cite{vanEkerenHeluani}.
    & \lean{VOA.VirasoroClassicalFreeness.classicallyFree}\\
    \addlinespace
    \textbf{Huang's tensor categories.} For strongly rational VOAs, the braided tensor category of modules; for
    $C_2$-cofinite VOAs, that of logarithmic modules. If the VOA is moreover pseudo-unitary: rigidity, the Verlinde
    formula and the congruence property of its modular representation.
    Source: \cite{HuangRigidity,HuangVerlinde,HuangLog,DongLinNg}.
    & \lean{VOA.Huang.TensorCat.HexagonClose.braidedCategory_stronglyRational_pbw},
    \lean{VOA.Huang.HuangRigid.rigid_ordinary},
    \lean{VOA.Classification.COne.huangVerlindeFusionInput_holds},
    \lean{VOA.Huang.CongruenceNamed.dlnCongruenceInput_holds}\\
    \addlinespace
    \textbf{Classification at $c=1$ (Gannon--Rayhaun).} The strongly rational, pseudo-unitary VOAs of central charge one
    are exactly the lattice VOAs
    $V_{2m}$, their orbifolds $V_{2m}^+$, and the three exceptional orbifolds $V_2^{A_4}$, $V_2^{S_4}$, $V_2^{A_5}$, each
    occurring once. Source: \cite{GannonRayhaun}.
    & \lean{VOA.Classification.COne.gannonRayhaun_classification}\\
    \addlinespace
    \textbf{Modular differential equations beyond $C_2$-cofiniteness.} $L_{-4/3}(\sltwo)$ is not $C_2$-cofinite, yet its
    characters converge, its torus conformal blocks form a finite-dimensional space, and its vacuum one-point functions
    satisfy a modular linear differential equation. The Arakawa--Kawasetsu theorem is proved for every VOA with
    finite-dimensional zeroth Poisson homology, the property of quasi-lisse VOAs that it uses. Source:
    \cite{ArakawaKawasetsu}.
    & \lean{VOA.Sl2NonLisse.not_C2Cofinite_simpleSl2QL},
    \lean{VOA.IsSqPrimary.exists_monicMLDEOn_all_analyticBlocks_of_etingofFinite_uncond}\\
    \midrule
    \multicolumn{2}{l}{\textbf{Monstrous moonshine}}\\
    \midrule
    \textbf{Conway--Norton for every automorphism of $V^\natural$ (Borcherds).} For every $g\in\mathrm{Aut}\,V^\natural$,
    $T_g$ is the normalised Hauptmodul of a genus-zero group containing some $\Gamma_0(N)$ and commensurable with
    $SL_2(\Z)$. No hypotheses. Qualitative form: $\Gamma_g$ is not matched with the Conway--Norton table. Source:
    \cite{ConwayNorton,Borcherds,CumminsGannon}.
    & \lean{VOA.Moonshine.ConwayNorton.conwayNorton}\\
    \addlinespace
    \textbf{Finiteness and rationality.} $\mathrm{Aut}\,V^\natural$ is finite, every automorphism has finite order,
    and every graded trace is rational. No hypotheses, no group theory.
    & \lean{VOA.Moonshine.ConwayNorton.M8.finite_autVnat}, \lean{VOA.Moonshine.ConwayNorton.M8.autVnat_isOfFinOrder},
    \lean{VOA.Moonshine.ConwayNorton.rationalTraces}\\
    \addlinespace
    \textbf{The Monster acts on $V^\natural$.} $V^\natural_2\cong\C\otimes\mathrm{SpaceB}$, AxiomMath's rational Griess
    algebra; every element of $\mathrm{Gmat}=\langle 2^{1+24}.\mathrm{Co}_1,\tau_B\rangle$ lifts to $V^\natural$ and
    satisfies moonshine; $\mathrm{Aut}\,V^\natural$ is isomorphic to the subgroup of $\mathrm{AutO}$ acting on
    $V^\natural$. No hypotheses beyond $g\in\mathrm{Gmat}$; Tits' theorem is not used. Source: \cite{FLM,Conway85,AxiomCFSG}.
    & \lean{VOA.Moonshine.ConwayNorton.conwayNorton_gmat},
    \lean{VOA.Moonshine.ConwayNorton.nonempty_autVnat_mulEquiv_liftable}\\
    \addlinespace
    \textbf{Monstrous moonshine for the Monster.} An isomorphism $\mathrm{AutO}\simeq\mathrm{Aut}\,V^\natural$
    compatible with the Griess algebras, and moonshine for every element of the Monster. No hypotheses in the
    package \texttt{MoonshineCFSG/}, which imports AxiomMath's generation theorem $\mathrm{Gmat}=\mathrm{AutO}$; in the
    library, given that theorem as the named hypothesis \lean{GmatEqAutO}.
    & \lean{VOA.Moonshine.ConwayNorton.monstrousMoonshine_unconditional} (package),
    \lean{VOA.Moonshine.ConwayNorton.monstrousMoonshine} (library)\\
    \midrule
    \multicolumn{2}{l}{\textbf{New results: rigidity}}\\
    \midrule
    \textbf{Virasoro minimal models.} Every graded first-order deformation of every minimal model $L(c_{p,q},0)$,
    $p,q\ge2$ coprime, is trivial. For the series $M(2,2n+1)$, rigidity was also proved by Butson and Nair by other
    methods \cite{ButsonNair}. The universal Virasoro VOA is proved not rigid, so the notion has content.
    & \lean{VOA.VirModule.minimal_graded_rigidity}, \lean{VOA.VirPoly.not_isRigidGraded_virasoroVOA}\\
    \addlinespace
    \textbf{Affine $\sltwo$ at admissible levels.} Every graded first-order deformation of $L_k(\sltwo)$ is trivial, at
    every admissible level. This contains Linshaw and Qi's results for $\sltwo$, at the positive integral levels and at
    $k=-4/3$ \cite{LinshawQi}.
    & \lean{VOA.AffineSl2GeneralRigidity.admissible_graded_rigidity}\\
    \midrule
    \multicolumn{2}{l}{\textbf{New results: classical freeness for the Deligne series}}\\
    \midrule
    \textbf{$A_2$, unconditional.} $L_{-3/2}(\slthree)$ is classically free: its $C_2$ algebra is the quotient by eight
    explicit quadrics, and the arc algebra of that quotient is isomorphic to the associated graded of the VOA.
    & \lean{VOA.DeligneAudit.DeligneSeries.A2_classicallyFree_final}\\
    \addlinespace
    \textbf{$D_4$, conditional on one published result.} The same for $L_{-2}(\mathfrak{so}_8)$, with $105$ quadrics,
    assuming a published dimension bound for certain fusion products of current-algebra modules \cite{KusLittelmann}.
    One step of its published proof needs repair; nothing suggests the bound is false.
    & \lean{VOA.DeligneAudit.DeligneSeries.D4_classicallyFree_of_inputT}\\
    \addlinespace
    \textbf{$G_2$, $F_4$, $E_6$, $E_7$, $E_8$, conditional on a conjecture.} The same at the other levels of the series,
    assuming the same published bound, a published result on Demazure modules \cite{Ravinder,ChariVenkatesh,KumarDCF,MathieuDCF}
    (proved in Lean for $D_4$, not yet for these types), and a character identity that we checked numerically to high
    order but that has no published proof.
    & \lean{VOA.DeligneAudit.DeligneSeries.G2_classicallyFree_final} and its $F_4$, $E_6$, $E_7$, $E_8$ analogues\\
    \midrule
    \multicolumn{2}{l}{\textbf{New results: graded unitarity}}\\
    \midrule
    \textbf{Universal affine $\sltwo$, unconditional.} For every real $k<-2$, $V^k(\sltwo)$ is graded unitary in the
    sense of \cite{ABLR}.
    & \lean{VOA.AffineSl2BelowCritical.universal_gradedUnitary_belowCritical_unconditional}\\
    \addlinespace
    \textbf{Transfer through the BRST coset.} If $L_k(\sltwo)$ at the boundary admissible level is graded unitary,
    then so is the minimal model $M(2,2n+1)$, for every $n\ge2$. This is a proved implication; its hypothesis is open.
    The proof formalizes the argument of Beem and Garner \cite{BG} on the actual BRST complex.
    & \lean{VOA.TwoSeries.JetPoincareTwo.twoSeries_gradedUnitary_of_affine_final}\\
  \end{longtable}
}

\paragraph{Also in the library.} Vertex operator superalgebras, with symplectic fermions and the $\beta\gamma$
system as examples, and the triplet algebra $W(2)$, proved simple, $C_2$-cofinite and not rational.

\section{Trust model}\label{sec:trust}

\paragraph{Base.} The library is built with Lean~v4.35.0-rc3 and Mathlib commit \texttt{b63f6e8a}; the package
\texttt{MoonshineCFSG/} with the same toolchain, the library and AxiomMath's library. Neither declares an axiom of its
own.
Production code contains no \lean{sorry}, \lean{admit}, \lean{native_decide}, \lean{stop} or
\lean{implemented_by}. We wrote an independent repository checker that rejects any of these, as well as any
\lean{axiom} declaration, and requires every library file to be reachable from the umbrella module. Large finite
verifications are done by kernel reduction (\lean{decide}, \texttt{decide +kernel}), never by the compiler. In addition,
we run \texttt{\#print axioms} on every headline declaration.

\paragraph{Definitions and theorems.} The real risk in a formalization of this kind is not an unsound proof but a true
theorem about the wrong object: a definition that is weaker than the textbook one, a structure that can only be
trivial, a conclusion satisfied by a degenerate witness.
The book is our main defense: it shows how each declaration is connected to the corresponding statement in the
literature.
For every theorem it prints the Lean declaration verbatim (a script checks that each
quoted block still occurs in the source at the recorded commit), decodes every convention and hypothesis in a
dictionary box, and states how the statement could have been wrong or vacuous and why it is not. Each chapter has
two referees: one for clarity, one for fidelity (statement against Lean, verbatim quotes). Before release, each
headline area received an independent \emph{release audit}: a fresh AI reviewer and human readers
compared the Lean statement with the source paper, traced every definition it depends on,
listed every explicit hypothesis, and re-printed the axioms. A final read of the book's release edition then brought every chapter into line with the audited wording.

\paragraph{The named hypothesis and its discharge.} In the library, \lean{GmatEqAutO} is a hypothesis of
\lean{monstrousMoonshine}, not an axiom: the library declares no axiom. Its content is the theorem
\lean{CFSG.Blueprint.NativeSporadic.MonsterModel.gmat_eq_griess} of the AxiomMath library \cite{AxiomCFSG}, whose proof
passes through the Stroth--Leon--Sims recognition of the Baby Monster. Our ported definitions
are textually AxiomMath's after renaming (the port-fidelity report \texttt{Notes/monster\_aut/CFSG\_PORT\_FIDELITY.md}
ships with the release). Textual identity was not enough, however. A Lean structure is a nominal type, so redeclared
copies of AxiomMath's structures are distinct types, and in v1.05 \lean{gmat_eq_griess} could not have proved
\lean{GmatEqAutO} even on a common toolchain. The package \texttt{MoonshineCFSG/} therefore uses AxiomMath's own structures (the
field $\mathbb F_4$, the spin coordinates, the Parker loop), its Golay code and its lattice map, so that every layer of our
\lean{Generated.Gmat} and \lean{AutO} is definitionally AxiomMath's. The kernel accepts
\lean{gmat_eq_griess} as a proof of \lean{GmatEqAutO}. Lean's elaborator times out unfolding the Griess product, so
the theorem \lean{gmatEqAutO} is added by a command that hands this one declaration directly to the kernel. Nothing
beyond the kernel is trusted, and \texttt{\#print axioms} reports it as for any theorem. The library keeps its own
copies of the structures, so that it does not depend on AxiomMath's library; the price is that its \lean{AutO} and the
package's are different, though textually identical, Lean objects, and the unconditional theorem is about the
package's, which is AxiomMath's own model. A script checks that the package's copies of library modules differ from
the originals only in the files that switch to AxiomMath's structures.

\paragraph{Verification of v1.1.} On Lean~v4.35.0-rc3, the library builds green (\texttt{lake build},
21{,}607 jobs on our development tree, which also contains the files held back from the release) with no
\lean{sorry}. Every released library file is byte-identical to a file of that build, except three umbrella files from
which only imports of held files were removed. The repository checker passes: every library
file is reachable from the umbrella module and there is no \lean{axiom} declaration. The axiom audit passes with 5{,}983
reports across 25 probes, and a separate script prints the explicit hypotheses and the axioms of 14 headline theorems,
each with the three standard axioms. The package \texttt{MoonshineCFSG/} with AxiomMath's library at commit
\texttt{5437cdc313} (whose definitions of \lean{Gmat} and \lean{AutO}, with their import closure, are byte-identical to
those of the commit we ported from) builds green: 31{,}417 jobs, including AxiomMath's proof of \lean{gmat_eq_griess}
built from source. \texttt{\#print axioms} for \lean{gmatEqAutO} and \lean{monstrousMoonshine_unconditional} gives
\lean{propext}, \lean{Classical.choice} and \lean{Quot.sound}; this also checks the axiom footprint of
\lean{gmat_eq_griess} itself. Every other result in Table~\ref{tab:scope}, including moonshine for every automorphism of
$V^\natural$ and for every element of $\mathrm{Gmat}$, is proved in the library and does not depend on the package.

\paragraph{The port.} Moving to Lean~v4.35 changed about $3{,}800$ files. The port adds no \lean{sorry}, \lean{axiom} or
other escape, and changes elaboration options on eight declarations only.
Changes in Mathlib forced new spellings of some statements, with the same mathematics (for instance, polynomial
coefficients are now written \lean{p.coeff m}, because Mathlib's monoid algebras became a structure), and three
definitions keep their type but are no longer definitionally equal to their former bodies. Every Lean excerpt in the
book was re-checked against the new sources.

\section{Size and structure}

The released library has 12{,}057 Lean files and 3{,}848{,}643 lines under \texttt{VOA/}. The package
\texttt{MoonshineCFSG/} has 86 further files: the bridge and the copies of the 85 library modules it needs.

\begin{center}\small
  \begin{tabular}{lrrl}
    \toprule
    Directory & Files & Lines & Content\\
    \midrule
    \texttt{Foundations} & 29 & 9{,}676 & formal calculus, PBW, Lie-theoretic input\\
    \texttt{Core} & 128 & 46{,}404 & VOA axioms, Borcherds identity, quotients\\
    \texttt{Modules} & 469 & 133{,}613 & weak/admissible/ordinary modules, Zhu algebra\\
    \texttt{Modularity} & 234 & 67{,}061 & modular forms, Zhu's theorem, $q$-traces\\
    \texttt{Examples} & 1{,}183 & 307{,}891 & constructions; rigidity; Arakawa--Dai; graded unitarity\\
    \texttt{AffineGeneral} & 354 & 21{,}656 & affine VOAs for every simple $\g$\\
    \texttt{LatticesAndCodes} & 26 & 8{,}551 & lattices, Golay code, Leech lattice\\
    \texttt{Moonshine} & 561 & 148{,}251 & $V^\natural$; monstrous moonshine; the Monster's action\\
    \texttt{Huang} & 931 & 184{,}162 & tensor categories, rigidity, Verlinde, congruence\\
    \texttt{Classification} & 1{,}756 & 1{,}628{,}058 & Gannon--Rayhaun $c=1$ classification\\
    \texttt{DeligneAudit} & 6{,}173 & 1{,}242{,}219 & Deligne series\\
    \texttt{Commutative}, \texttt{Filtered}, \texttt{Geometric} & 201 & 45{,}649 & commutative algebra, Li filtration\\
    umbrella files & 12 & 5{,}452 & \\
    \bottomrule
  \end{tabular}
\end{center}

A large part of the line count is certificate data. The 167 files larger than 200\,KB contain 1{,}518{,}787 lines
(about 40\% of the total), mostly exact rational data for the $c=1$ orbifold computations; they are checked by kernel
reduction (about 39{,}600 \texttt{decide +kernel} calls), and the book describes them collectively --- what is checked,
how, and at what size --- while presenting the theorems that consume them individually.

The release edition of the book has 2{,}807 pages, and its brief edition 258 pages; every Lean block quoted in
the brief is checked verbatim against the source. Part~VIII of the book, on monstrous moonshine, has seven chapters,
each refereed twice.

\section{Method}

The formalization was human-directed and carried out with the assistance of AI coding agents (Anthropic's
Claude Opus and Claude Fable, and OpenAI's Codex), while Python scripts generated the certificate data that the Lean
kernel re-checks and Python and Mathematica supplied independent numerical cross-checks. The authors chose the
targets, the definitions and conventions, and the release scope; agents wrote all of the Lean code and the first
drafts of the book chapters. Three practices shaped the result. First, referee loops: every chapter and every
headline went through independent review by several open-weight and proprietary AI models, and statements were
corrected or narrowed when a review found an over-claim. Second, a strict reporting rule: a result counts only when the aggregate build is green, and
``unconditional'' is reserved for statements whose hypothesis closure is fully discharged. Third, agents were required to record every dead
end, which let later agents narrow their search for proofs:
routes proved false, vacuous or circular are kept, with their counterexamples, in a project log, and some are shown as
examples in the book's verification chapter.

\section{In preparation}\label{sec:prep}

\begin{itemize}
  \item \textbf{The table-level Conway--Norton statement}, identifying each $\Gamma_g$ with the group in Conway and
    Norton's table \cite{ConwayNorton}.
  \item \textbf{Generalized moonshine} (Norton), for the twisted modules of $V^\natural$.
  \item \textbf{Deligne series beyond $A_2$.} Removing the assumptions of the conditional theorems; for $D_4$ the
    book (Chapter~L27) describes the current state. Newer Lean files on the Deligne series are held back from v1.1
    and will be released when classical freeness for the whole series is unconditional.
  \item \textbf{Graded unitarity in the other direction}, from the minimal models to the affine algebra, through
    Adamovi\'c's free-field realization (inverse Drinfeld--Sokolov reduction).
  \item \textbf{Affine $\slthree$}, towards its rigidity.
  \item \textbf{Removing pseudo-unitarity.} Proving rigidity, the Verlinde formula and the congruence property for
    every strongly rational VOA, as in \cite{HuangRigidity,HuangVerlinde,DongLinNg}, without assuming
    pseudo-unitarity.
\end{itemize}

\section*{Acknowledgements}
We thank Brandon Rayhaun, Avia Raviv-Moshe and Chris Beem for useful discussions. We also thank the Lean and
Mathlib communities. The work of LR is supported in part by NSF grant PHY-2513893 and by the Simons Foundation
grant 681267 (Simons Investigator Award).


\end{document}